\documentclass[11pt]{article}
\usepackage{amsfonts}
\usepackage{amsmath}
\usepackage{amssymb}
\usepackage{latexsym,color}
\usepackage{amsbsy}
\usepackage{graphicx}

\newtheorem{thm}{Theorem}[section]
\newtheorem{con}{Conjecture}[section]
\newtheorem{lemma}[thm]{Lemma}
\newtheorem{cor}[thm]{Corollary}
\newtheorem{pro}[thm]{Proposition}
\newtheorem{example}[thm]{Example}
\newtheorem{definition}[thm]{Definition}

\newtheorem{remark}[thm]{Remark}
\newtheorem{Algorithm}[thm]{Algorithm}

\newcommand{\comment}[1]{}

\newcommand{\N}{{\mathbb N}}
\newcommand{\Bq}{B_q}
\newcommand{\Cq}{C_q}
\newcommand{\Bfq}{B_{{\mathbb F}_q}}
\newcommand{\Cfq}{C_{{\mathbb F}_q}}
\newcommand{\Fq}{{\mathbb F}_q}
\newcommand{\C}{{\mathbb C}}

\newcommand{\pr}{\parallel}

\newcommand{\kq}[1]{{{ [#1]_q }}}
\newcommand{\qb}[2]{{{ {#1}\brack {#2}}}_q }

\newcommand{\ncom}{\newcommand}
\ncom{\ns}{\normalsize}
\ncom{\la}{\lambda}
\ncom{\bm}{\boldmath}
\ncom{\noi}{\noindent}
\ncom{\bq}{\begin{equation}}
\ncom{\eq}{\end{equation}}
\ncom{\beqn}{\begin{eqnarray*}}
\ncom{\eeqn}{\end{eqnarray*}}
\ncom{\ba}{\begin{array}}
\ncom{\ul}[1]{{#1}^-}

\ncom{\ea}{\end{array}}
\ncom{\beq}{\begin{eqnarray}}
\ncom{\eeq}{\end{eqnarray}}
\ncom{\nno}{\nonumber}
\ncom{\hs}{\mbox{\hspace{.25cm}}}
\ncom{\rar}{\rightarrow}
\ncom{\Rar}{\Rightarrow}
\ncom{\noin}{\noindent}
\ncom{\bc}{\begin{center}}
\ncom{\ec}{\end{center}}
\ncom{\sz}{\scriptsize}
\ncom{\fpd}{\Phi(\pi^{'})}
\ncom{\fp}{\Phi(\pi) }
\ncom{\nk}{\left< \begin{array}{c}
                       n\\k \end{array} \right>}
\ncom{\nd}{1^{'},2^{'},\cdots,n^{'}}
\ncom{\R}{I\!\!R}
\ncom{\de}{\bigtriangleup (F_{2n},\leq)}
\ncom{\del}{\bigtriangleup}
\ncom{\cov}{<\!\!\!\!\cdot }
\ncom{\bt}{\begin{thm}}
\ncom{\bcon}{\begin{con}}
\ncom{\et}{\end{thm}}
\ncom{\econ}{\end{con}}
\ncom{\bl}{\begin{lemma}}
\ncom{\el}{\end{lemma}}
\ncom{\bco}{\begin{cor}}
\ncom{\ds}{\displaystyle}
\ncom{\eco}{\end{cor}}
\ncom{\bp}{\begin{pro}}
\ncom{\ep}{\end{pro}}
\ncom{\bex}{\begin{example}}
\ncom{\eex}{\end{example}}
\ncom{\bd}{\begin{definition}}
\ncom{\ed}{\end{definition}}
\ncom{\brm}{\begin{remark}}
\ncom{\erm}{\end{remark}}
\ncom{\bal}{\begin{Algorithm}}
\ncom{\eal}{\end{Algorithm}}
\ncom{\ol}{\overline}

\ncom{\pf}{\noi {\bf Proof  }}
\ncom{\be}{\begin{enumerate}}
\ncom{\ee}{\end{enumerate}}
\ncom{\s}{\subset}
\ncom{\T}{{\cal T}}
\ncom{\B}{{\cal B}}
\ncom{\A}{{\cal A}}

\title{\Large{{\bf 
Counting spanning trees of the hypercube and its $q$-analogs by explicit
block diagonalization
}}}

\author{{ { Murali K. Srinivasan}} \\
{\em  \normalsize{Department of Mathematics}}\\
{\em  \normalsize{Indian Institute of Technology, Bombay}}\\
{\em  \normalsize{Powai, Mumbai 400076, INDIA}}\\
{\bf  \texttt{mks@math.iitb.ac.in}}}

\begin{document}
\date{}
\maketitle

\begin{abstract} The number of spanning
trees of a graph $G$ is  called the {\em complexity} of $G$ and is denoted
$c(G)$. Let $C(n)$ denote the {\em (binary) hypercube} of dimension $n$. A
classical result in enumerative combinatorics (based on explicit
diagonalization) states that $c(C(n)) = \prod_{k=2}^n (2k)^{n\choose k}$. 

In this paper we use the explicit block diagonalization methodology to
derive formulas for the complexity of two $q$-analogs of $C(n)$, the {\em
nonbinary hypercube} $\Cq(n)$, defined for $q\geq 2$, and the {\em vector
space analog of the hypercube} $\Cfq(n)$, defined for prime powers $q$.

We consider the nonbinary and vector space analogs of the
Boolean algebra. We show the existence, in both cases, of a graded Jordan
basis (with respect to the up operator) that is orthogonal
(with respect to the standard inner product) and we write down
explicit formulas for the ratio of the lengths of the successive vectors in
the Jordan chains (i.e., the singular values). With respect to (the
normalizations of) these bases
the Laplacians of $\Cq(n)$ and $\Cfq(n)$ block diagonalize, with
quadratically many distinct blocks in the nonbinary case and linearly many
distinct blocks in the vector space case, and
with each block
an explicitly written down real, symmetric, tridiagonal matrix of known
multiplicity and size at most $n+1$. In the nonbinary case we further
determine the eigenvalues of the blocks, by explicitly writing out the
eigenvectors, yielding an explicit formula for $c(\Cq(n))$ (this proof
yields new information even in the binary case). In the vector
space case we have been unable to determine the eigenvalues of the blocks
but we give a useful formula for $c(\Cfq(n))$ involving ``small" determinants (of
size at most $n$).

\end{abstract}

{{\bf  \section {  Introduction }}}  

Explicit block diagonalization was pioneered in the classic paper of
Schrijver {\bf\cite{s}} to improve the polynomial time computable 
Delsarte linear programming bound
on binary code size by using semidefinite programming. 
In this paper we apply this methodology to study 
two counting problems. 

Suppose we have a family $\{M(n)\}_{n\geq 1}$ of real, symmetric matrices,
where the size $s(n)$ of $M(n)$ is exponential in $n$. We are interested in
a formula for $D(n) = \,\mbox{det}\,M(n)$. In many combinatorial situations
the eigenvalues have large multiplicity due to the presence of symmetry. 
Suppose we find that:

(i) $M(n)$ has $p(n)$ distinct eigenvalues, where $p(n)$ is bounded by 
a polynomial in $n$.

(ii) We can determine the eigenvalues $\lambda_n(1),\ldots
,\lambda_n(p(n))$ of $M(n)$.

(iii) We can determine the multiplicity $m_n(i)$ of the
eigenvalue $\lambda_n(i),\;i=1,\ldots ,p(n)$ of $M(n)$.

Under these conditions it is clear that
\beq \label{ed} D(n)&=& \prod_{i=1}^{p(n)} \lambda_n(i)^{m_n(i)}\eeq
is a satisfactory formula for $D(n)$ 
and we say that (\ref{ed})
has been obtained by {\em explicit diagonalization}. 

Explicit diagonalization is the best case of   
explicit block diagonalization.
Suppose we are able to find a basis under which:

(a) $M(n)$ is in block diagonal form with possibly repeated blocks, but
the number $p(n)$ of distinct blocks is bounded by a polynomial in $n$.

(b) We can explicitly write down the distinct blocks
$B_n(1),\ldots ,B_n(p(n))$ of $M(n)$.

(c) The size $s_n(i)$ of the block $B_n(i)$ is bounded by 
a polynomial in $n$, for all $i=1,\ldots ,p(n)$.

(d) We can determine the multiplicity $m_n(i)$ of the
block $B_n(i)$ of $M(n)$, $i=1,\ldots ,p(n)$.

Under these conditions we have
\beq \label{ebd} D(n)&=& \prod_{i=1}^{p(n)} \,(\mbox{det}\,B_n(i))^{m_n(i)}\eeq
and we say that (\ref{ebd})
has been obtained by {\em explicit block diagonalization}. Note that 
conditions (a) and (c) taken together imply condition (i) 
in the paragraph above.
In this sense (\ref{ebd}) is much closer to (\ref{ed})
than to the formula $D(n) = \,\mbox{det}\,M(n)$, although
it is not as good as (\ref{ed}). 
Also note that if we can now
explicitly  determine the eigenvalues of the blocks $B_n(i)$ then we would
have achieved explicit diagonalization.

In this paper we give formulas of type (\ref{ebd}) above for two counting
problems. For one of these problems we can actually determine the
eigenvalues of the blocks, yielding a formula of type (\ref{ed}) above. We
now state our results.

The number of spanning trees of a graph $G$ is called the {\em
complexity} of $G$ and denoted $c(G)$. The {\em (binary) hypercube} $C(n)$ 
is the regular graph of degree $n$ whose vertex set is
the set of all $2^n$ subsets of the $n$-set $[n]=\{1,2,\ldots ,n\}$ and 
where two
subsets $X,Y\subseteq [n]$ are connected by an edge iff 
$X\subseteq Y$ or $Y\subseteq X$, and $||X|-|Y|| = 1$.
A beautiful classical result in enumerative combinatorics (based on explicit
diagonalization, see Example
5.6.10 in {\bf\cite{s3}})  states that
\beq \label{hc}
c(C(n)) &=& \prod_{k=2}^{n} (2k)^{n\choose k}
\eeq
%

We use the explicit block diagonalization methodology to
derive formulas for the complexity of two nonregular
$q$-analogs of $C(n)$, the {\em nonbinary
hypercube} $\Cq(n)$, defined for $q\geq 2$, and the {\em vector space analog of
the hypercube} $\Cfq(n)$, defined for prime powers $q$. 

The starting point of this paper is an alternative formulation and
interpretation of formula (\ref{hc}) for $c(C(n))$. 
There are two primary reasons for this.
Firstly,
we do not know a vector space  analog of (\ref{hc}) and we need the present 
approach to develop 
a unified theory covering all three cases $C(n), \Cq(n)$, and
$\Cfq(n)$.
Secondly, even in the nonbinary case, where an
analog of (\ref{hc}) is available (based on the product structure of
$\Cq(n)$), 
our approach yields more information. For instance,
one of our main results, Theorem \ref{mt4} in Section 3, is a
natural byproduct of the present approach (see ``Background and
Motivation" at the end of this section).

We have 
\beq \nonumber 
c(C(n)) &= & \frac{1}{2^n} \left\{\ds {\prod_{k=1}^n (2k)^{n\choose k}
}\right\} \\ 
\nonumber &=& \frac{1}{2^n} \left\{ {\ds{\prod_{k=1}^{n} (2k) }}\right\} 
\left\{ 
{\ds{\prod_{k=1}^{\lfloor n/2 \rfloor}}} 
\left({\ds{\prod_{j=k}^{n-k}}} 
(2j)\right)^{{n\choose
k}-{n\choose{k-1}} } \right\}\\
   \label{hcmks}    &=&   n! 
\left\{ 
{\ds{\prod_{k=1}^{\lfloor n/2 \rfloor}}} 
\left({\ds{\prod_{j=k}^{n-k}}} 
2j)\right)^{{n\choose
k}-{n\choose{k-1}} } \right\}
\eeq
To see the equivalence of the first and second lines above note that,
for $1\leq j\leq n/2$, the exponent of $2j$ in the numerator of the first line is
$n\choose j$ and in the numerator of the second line is also ${n\choose j}=1 + {n\choose 1}- {n\choose
0}+\cdots +{n\choose j}-{n\choose {j-1}}$. Since ${n\choose k} = {n\choose
{n-k}}$ the same conclusion holds for $n/2 \leq j \leq n$. The $q=2$ case of
the proof of Theorem \ref{stnb} below provides a linear algebraic 
interpretation to (\ref{hcmks}).

Let $q\geq 2$. Define $\Bq(n)$ to be the set of all pairs $(X,f)$, where
$X\subseteq [n]$ and $f:X\rar \{1,2,\ldots ,q-1\}$. The {\em nonbinary
hypercube} $\Cq(n)$ is the graph whose vertex set is $\Bq(n)$ and where
two vertices $(X,f)$ and $(Y,g)$ are connected by an edge iff 
$X\subseteq Y$ or $Y\subseteq X$, $||X|-|Y|| = 1$, and $f,g$ agree on $X\cap
Y$.
An equivalent way to define $\Cq(n)$ is as follows:
the vertex set is the set of all $n$-tuples $a=(a_1,\ldots ,a_n)$, where $a_i \in
\{0,1,\ldots ,q-1\}$ for all $i$. Define $\mbox{supp}(a)=\{i: a_i \not=
0\}$ and connect $a$ and $b=(b_1,\ldots ,b_n)$ by an edge iff 
$\mbox{supp}(a) \subseteq \mbox{supp}(b)$ or $\mbox{supp}(b) \subseteq
\mbox{supp}(a)$, $||\mbox{supp}(a)|-|\mbox{supp}(b)|| = 1$, and 
$a_i = b_i$ for $i\in \mbox{supp}(a) \cap \mbox{supp}(b)$. 
In this description it is clear that $|\Bq(n)|= q^n$ (we use this
description of $\Bq(n)$ in Section 5). 
Note that $\Cq(n)$ is different from what is usually called the nonbinary
Hamming graph $H_q(n)$. Both 
$H_q(n)$ and $\Cq(n)$ have the same vertex set
(consisting of all $n$-tuples with entries in $\{0,1,\ldots ,q-1\}$) but 
two vertices are connected by an edge in $H_q(n)$ iff they differ in exactly
one coordinate whereas they are connected by an edge in $\Cq(n)$ iff they
differ in exactly one coordinate and one of the $n$-tuples is zero in that
coordinate.  
Note also that 
$\Cq(n)$ is nonregular. However, the degree $|X| + (q-1)(n- |X|)$ 
of the vertex $(X,f)$ in $\Cq(n)$ depends only on $|X|$.
This property will prove useful in counting the spanning trees of $\Cq(n)$.

In Section 3 we prove
the following explicit formula for $c(\Cq(n))$.
For $0\leq k \leq n$, define $\ul{k} = \mbox{ max}\{0, 2k -  n\}$. 
Clearly $0\leq \ul{k} \leq k$ and $k\leq n + \ul{k} - k$.  
\bt \label{stnb} We have
$$
c(\Cq(n)) = n! \left\{
{\ds{\prod_{k=1}^{n}}}\,
{\ds{\prod_{l= \ul{k}}^{k}}}
\left({\ds{\prod_{j=k}^{n+l-k}}}
(qj - (q-1)l)\right)^{(q-2)^l{n\choose l}\left\{{{n-l}\choose
{k-l}}-{{n-l}\choose{k-l-1}}\right\}} \right\}
$$
\et
Note that when $q=2$ only terms with $l=0$ survive (since $0^0 = 1$) and 
then the formula above reduces to $c(C(n))$.

Now let $q$ be a prime power. Let $\Bfq(n)$ denote 
the set of all subspaces of an $n$-dimensional vector space over the finite
field $\Fq$ and set $G_n(q)=|\Bfq(n)|$. The {\em Galois numbers} $G_n(q)$
satisfy the recursion $G_{n+1}(q) = 2G_n(q) + (q^n - 1)G_{n-1}(q),\;n\geq 1$
(see Exercise 1.73 in {\bf\cite{a}}). The {\em $q$-binomial or Gaussian coefficient}
$\qb{n}{k}$ denotes the number of $k$-dimensional subspaces in $\Bfq(n)$.

The {\em vector space analog $\Cfq(n)$ of the hypercube} 
is the graph whose vertex set is $\Bfq(n)$,
and where  
subspaces $X,Y\in \Bfq(n)$ are connected by an edge iff 
$X\subseteq Y$ or $Y\subseteq X$, and $|\mbox{dim}(X) - \mbox{ dim}(Y)| = 1$.
Note that the graph
$\Cfq(n)$ is nonregular. However, the degree $\kq{k}+ \kq{n-k}$ (where, for 
$k\in \N \;(= \{0,1,2,\ldots \})$, we set 
$\kq{k} =1+q+q^2+\cdots +q^{k-1}$) 
of the vertex $X$ in $\Cfq(n)$ depends only on $k=\mbox{ dim}(X)$.
This property will prove useful in counting the spanning trees of $\Cfq(n)$.

In Section 4 we prove the following formula for $c(\Cfq(n))$. It is similar to 
formula (\ref{hcmks}) for $c(C(n))$, with the important
difference that the explicit term 
$\prod_{j=k}^{n-k}(2j)$ is replaced by a recursive calculation.
It would be best if this recurrence were replaced by an explicit term,
yielding a ``closed form"  formula for $c(\Cfq(n))$  
(this would amount to
explicitly determining the eigenvalues of the Laplacian of $\Cfq(n)$, see
Section 4). 
We do not
know how to do this. 

Let $k,n \in \N$ with $k\leq n/2$. For
$k \leq j \leq n-k+1$, define polynomials $F_q(n,k,j)$ in $q$, having
integral coefficients, using the
following recursion: 
\beqn
F_q(n,k,n-k+1)&=&1 \\
F_q(n,k,n-k) &=& \kq{k} + \kq{n-k}
\eeqn
and, for $k\leq j < n-k$,
\beq \label{dr}
\nonumber \lefteqn{F_q(n,k,j) =} \\  
&(\kq{j} + \kq{n-j})F_q(n,k,j+1) -
(q^k\kq{j+1-k}\kq{n-k-j})F_q(n,k,j+2).&
\eeq

\bt \label{stvs} We have
$$c(\Cfq(n)) = \kq{1}\kq{2}\cdots \kq{n} \left\{
{\ds{\prod_{k=1}^{\lfloor n/2 \rfloor}}}
F_q(n,k,k)^{ \qb{n}{k}  - \qb{n}{k-1} } \right\}$$
\et
The following table gives the first five values of $c(\Cfq(n))$.
\beqn
c(\Cfq(1)) &= & 1\\
c(\Cfq(2)) &= & \kq{2} 2^q \\
c(\Cfq(3)) &= & \kq{2} \kq{3} (4+3q+q^2)^{q(1+q)}\\
c(\Cfq(4)) &= & \kq{2} \kq{3} \kq{4}
          (8+12q+12q^2+10q^3+4q^4+2q^5)^{q(1+q+q^2)}\\
          &&\times (2+2q)^{q^2(q^2+1)}\\
c(\Cfq(5)) &= &
\kq{2} \kq{3} \kq{4} \kq{5}
F_q(5,1,1)^{q(1+q)(1+q^2)}\\
&& \times F_q(5,2,2)^{q^2(1+q+q^2+q^3+q^4)}
\eeqn
where
$F_q(5,2,2) = 4+8q+7q^2+4q^3+q^4$ 
and 
$F_q(5,1,1) =$
$$16+36q+53q^2+65q^3+69q^4+58q^5+42q^6+26q^7+13q^8+5q^9+q^{10}.$$

Since the degree (in $q$) of the polynomials $\kq{j} + \kq{n-j}$ and
$q^k\kq{j+1-k}\kq{n-k-j}$ are both $\leq n-1$, it follows 
by an easy induction, using (\ref{dr}), that the 
degree of $F_q(n,k,j)$ is $\leq (n-1)(n-k+1-j)$ 
(an exact formula for the degree is given in Section 4).
Thus all the
polynomials $F_q(n,k,j)$ can be computed efficiently (in time polynomial
in $n$).

{\bf Remark} A special case of a result of Butler {\bf\cite{b}} shows that, as 
a polynomial in
$q$, $\qb{n}{k}  - \qb{n}{k-1},\;k\leq n/2$ has nonnegative integral
coefficients. Data suggest that the polynomials $F_q(n,k,j)$ have
nonnegative coefficients and are unimodal, for all $n,k,j$. 
We do not study this problem in the present paper. 
 
\noi
{\bf Background and motivation} The present approach has two main steps.

(i) {\em Explicit block diagonalization of the Laplacian:} 
The Laplacian matrices of $\Cq(n)$ and $\Cfq(n)$ are of 
exponential sizes $q^n$ and $G_n(q)$ (note that the size
depends on both $n$ and $q$). In Section 2 
we interpret the graphs $C(n), \Cq(n), \mbox{ and }\Cfq(n)$ as  the Hasse
diagrams of three natural graded posets, namely, the 
Boolean algebra, the nonbinary analog of the Boolean algebra, and the vector
space analog of the Boolean algebra.
We summarize (without proofs) results on the up operator on 
these posets (the results are proved in Section 5).
These results give unitary matrices, of respective sizes $2^n, q^n$, 
and $G_n(q)$, conjugating by which block diagonalizes the Laplacians of
$C(n), \Cq(n),$ and $\Cfq(n)$, with
polynomially many distinct blocks (in fact quadratic in the $\Cq(n)$ case
and linear in the other two cases),
where the multiplicity of each block is known and where
each block is an explicitly written down real, symmetric, tridiagonal matrix of
size at most $n+1$ (and independent of $q$). 
Since only the entries of the
blocks, and not their sizes, depend on $q$ we can now treat 
$q$ symbolically. 

The main inspiration and motivation for the results in Section 2 
are the papers of Schrijver {\bf\cite{s}} and Gijswijt, Schrijver, and
Tanaka {\bf\cite{gst}}.
Schrijver's explicit block diagonalization of the
commutant of the symmetric group action on the Boolean algebra (=
Terwilliger algebra of the binary Hamming scheme)
was derived from
Theorem \ref{mt1} below in {\bf\cite{s1}}. Similarly, the explicit block
diagonalization of the Terwilliger algebra of the nonbinary Hamming scheme
worked out in {\bf\cite{gst}} can be derived from
Theorem \ref{mt2}  and Theorem \ref{mt3} can be used to explicitly block
diagonalize the commutant of the $GL(n,\Fq)$ action on $V(\Bfq(n))$. 
Likewise, Theorems \ref{mt1}, \ref{mt2}, and
\ref{mt3} also provide a unified approach to the explicit diagonalization of
the Bose-Mesner algebras of the (binary) Johnson scheme and its nonbinary
and vector space analogs {\bf\cite{bi,tag}}. 
We do not discuss this topic in the present paper (we hope to write this
down in an expository paper {\bf\cite{s2}}).

(ii) {\em Explicit diagonalization of the blocks:} Using the matrix tree
theorem step (i) above yields a formula of type (\ref{ebd}) for the complexity.
To obtain an explicit formula
we need to
determine the eigenvalues of each distinct block appearing in the block
diagonalization in step (i). 
We carry out this step for $\Cq(n)$ in Section 3 by explicitly writing out
the eigenvectors of the blocks, 
yielding a proof of Theorem \ref{stnb}.
In the binary case, the eigenvectors of the Laplacian produced
by this proof are different from that given in the standard proof (Example
5.6.10 in {\bf\cite{s3}}). This is easily seen from the fact that all the
eigenvectors in the proof in {\bf\cite{s3}} have support of cardinality
$2^n$ (i.e., have nonzero components in each of the standard coordinates)
whereas that is not the case here. 

We have been unable to
carry out this step for $\Cfq(n)$
and this accounts for the nonexplicit nature of the formula in
Theorem \ref{stvs} (proved in Section 4). 
The blocks being symmetric,
tridiagonal their determinants can be easily calculated recursively,
yielding the polynomials $F_q(n,k,j)$.

{\textcolor{black} {\bf \section {Orthogonal Jordan chains and singular
values}}}

 A (finite) {\em graded poset} is a (finite) poset $P$ together with a
{\em rank function}
$r: P\rar \N$ such that if $p'$ covers $p$ in $P$ then
$r(p')=r(p)+1$. The {\em rank} of $P$ is $r(P)=\mbox{max}\{r(p): p\in P\}$
and,
for $i=0,1,\ldots ,r(P)$, $P_i$ denotes the set of elements of $P$ of rank
$i$. For a subset $S\subseteq P$, we set $\mbox{rankset}(S) = \{r(p):p\in
S\}$.

For a finite set $S$, let $V(S)$ denote the complex vector space with $S$ as
basis. Let $P$ be a graded poset with $n=r(P)$. Then we have
$ V(P)=V(P_0)\oplus V(P_1) \oplus \cdots \oplus V(P_n)$ (vector space direct
sum).
An element $v\in V(P)$ is {\em homogeneous} if $v\in V(P_i)$ for some $i$,
and we extend the notion of rank to homogeneous elements by writing
$r(v)=i$. Given an element $v\in V(P)$, write $v=v_0 + \cdots +v_n,\;v_i \in
V(P_i),\;0\leq i \leq n$. We refer to the $v_i$ as the {\em homogeneous
components} of $v$. A subspace $W\subseteq V(P)$ is {\em homogeneous} if it
contains the homogeneous components of each of its elements. For a
homogeneous subspace $W\subseteq V(P)$ we set $\mbox{rankset}(W)=\{r(v) : v
\mbox{ is a homogeneous element of } W\}$.

The {\em
up operator}  $U:V(P)\rar V(P)$ is defined, for $p\in P$, by
$U(p)= \sum_{p'} p'$,
where the sum is over all $p'$ covering $p$.
A {\em graded Jordan chain} in $V(P)$ is a sequence 
\beq \label{gjc}
&s=(v_1,\ldots ,v_h)&
\eeq 
of nonzero homogeneous elements of $V(P)$
such that $U(v_{i-1})=v_i$, for
$i=2,\ldots h$, and $U(v_h)=0$ (note that the
elements of this sequence are linearly independent, being nonzero and of
different ranks). We say that $s$ {\em
starts} at rank $r(v_1)$ and {\em ends} at rank $r(v_h)$. 
A {\em graded Jordan basis} of $V(P)$ is a basis $V(P)$
consisting of a disjoint union of graded Jordan chains 
in $V(P)$.

The graded Jordan
chain (\ref{gjc}) is said to be a {\em symmetric Jordan chain} (SJC) if
the sum of the starting and ending ranks of $s$ equals $r(P)$, i.e., 
$r(v_1) + r(v_h) = r(P)$ 
if $h\geq
2$, or $2r(v_1)= r(P)$ if $h=1$.
A {\em symmetric Jordan basis} (SJB) of $V(P)$ is a basis of $V(P)$
consisting of a disjoint union of symmetric Jordan chains 
in $V(P)$. 

The graded Jordan
chain (\ref{gjc}) is said to be a {\em semisymmetric Jordan chain} (SSJC) if
the sum of the starting and ending ranks of $s$ is $\geq r(P)$.
A {\em semisymmetric Jordan basis} (SSJB) of $V(P)$ is a basis of $V(P)$
consisting of a disjoint union of semisymmetric Jordan chains 
in $V(P)$. 
An SSJB is said to be
{\em rank complete} if it contains graded Jordan chains starting at rank $i$
and ending at rank $j$, for all $0\leq i\leq j \leq r(P),\;i+j\geq r(P)$.

Let $\langle , \rangle$ denote the standard inner product on $V(P)$,
i.e.,
$\langle p,p' \rangle =
\delta (p,p')$ (Kronecker delta), for $p,p'\in P$.
The {\em length} $\sqrt{\langle v, v \rangle }$ of $v\in V(P)$ is denoted
$\pr v \pr$.

Suppose we have an orthogonal graded Jordan basis $J(n)$ of $V(P)$.
Normalize the vectors in $J(n)$ to get an
orthonormal basis $J'(n)$.
Let $(v_1,\ldots ,v_h)$ be a graded Jordan chain in $J(n)$.
Put
$v_u' = \frac{v_u}{\pr v_u \pr}$ and $\alpha_u = \frac{\pr v_{u+1} \pr}{\pr
v_{u} \pr},\;1\leq u \leq h$ (we take $v_0'=v_{h+1}'=0$). 

We have,
for $1\leq u \leq h$,
\beq \label{trick}
U(v_{u}')=\frac{U(v_{u})}{\pr v_{u} \pr}=\frac{v_{u+1}}{\pr v_{u}
\pr}=\alpha_{u} v_{u+1}'.&&
\eeq
Thus the matrix of $U$ with respect to (wrt) $J'(n)$ is in block diagonal form, with a block
corresponding to each (normalized) graded Jordan chain in $J(n)$, and with
the block corresponding to $(v_1',\ldots ,v_h')$ above being a lower triangular
matrix with subdiagonal $(\alpha_1 ,\ldots ,\alpha_{h-1})$ and $0$'s
elsewhere.

The {\em
down operator}  $D:V(P)\rar V(P)$ is defined, for $p\in P$, by
$D(p)= \sum_{p'} p'$,
where the sum is over all $p'$ covered by  $p$.
Note that the matrices, in the standard basis,
of $U$ and $D$ are real and transposes of
each other. Since $J'(n)$ is orthonormal
wrt the standard inner product, it follows that the matrices of
$U$ and $D$, in the basis $J'(n)$, must be adjoints of each other.
Thus, for $0\leq u \leq h-1$,  we must have (using (\ref{trick}) and the
previous paragraph),
\beq \label{trick1}
D(v_{u+1}')=\alpha_{u} v_{u}'.&&
\eeq
In particular, the subspace spanned by $\{v_1,\ldots
,v_h\}$ is closed under $U$ and $D$.

Another useful observation is the following: take scalars $\beta_0 ,\beta_1
,\ldots ,\beta_{r(P)}$ and define the operator ${\cal B}: V(P)\rar V(P)$ by
${\cal B}(p) = \beta_{r(p)} p,\;p\in P$. Since each element of the graded
Jordan chain $(v_1,\ldots ,v_h)$ is homogeneous, it follows from the
definition of ${\cal B}$ that the subspace spanned by $\{v_1 ,\ldots ,v_h\}$
is closed under $U,D$ and ${\cal B}$.

In this paper we consider three graded posets. 
The {\em Boolean algebra} $B(n)$ is the graded poset of subsets (under
inclusion) of
$[n]=\{1,2,\ldots ,n\}$, with rank of a subset given by cardinality. 

\bt \label{mt1}
There exists a SJB $J(n)$ of $V(B(n))$ such that 

\noi (i) The elements of $J(n)$ are orthogonal with respect to
$\langle , \rangle$ (the standard inner product).

\noi (ii) {\em (Singular Values)} Let $0\leq k \leq n/2 $ and let
$(x_k,\ldots ,x_{n-k})$ be any SJC
in $J(n)$ starting at rank $k$ and ending at rank $n-k$. Then
we have, for $k\leq u < n-k$,
\beq \label{mks1}
\frac{\pr x_{u+1} \pr}{\pr x_u \pr} & = &
\sqrt{(u+1-k)(n-k-u)}
\eeq
\et

We now consider two $q$-analogs of the Boolean algebra, the nonbinary analog
and the vector space analog.

Partially order $\Bq(n)$ as follows:
$(X,f)\leq (Y,g)$ provided $X\subseteq Y$ and $f(a)=g(a),\;a\in X$. It is
easy to see that $\Bq(n)$ is a rank-$n$ graded poset with rank of $(X,f)$
given by cardinality of $X$.
We can also think of the
poset $\Bq(n)$ as the product $N_q \times \cdots \times N_q$ ($n$ factors),
where $N_q$ is the poset on $\{0,1,\ldots ,q-1\}$ with the covering relations $0 <
i$, $i\in \{1,\ldots ,q-1\}$.

For a SSJC $s$ in $V(\Bq(n))$, starting at rank $i$ and ending at rank $j$,
we define the {\em offset} of $s$ to be $i+j-n$. It is easy to see that if
an SSJC starts at rank $k$ then its offset $l$ satisfies  
$\ul{k}\leq l \leq k$ and the chain ends at rank $n+l-k$.

\bt \label{mt2}
There exists a rank complete SSJB $J(n)$ of $V(\Bq(n))$ such that 

\noi (i) The elements of $J(n)$ are orthogonal with respect to
$\langle , \rangle$ (the standard inner product).

\noi (ii) {\em (Singular Values)} Let $0\leq k \leq n$, $\ul{k}\leq l \leq k$ 
and let
$(x_k,\ldots ,x_{n+l-k})$  be any SSJC
in $J(n)$ starting at rank $k$ and having offset $l$. Then
we have, for $k\leq u < n+l-k$,
\beq \label{mks2}
\frac{\pr x_{u+1} \pr}{\pr x_u \pr} & = &
\sqrt{(q-1)(u+1-k)(n+l-k-u)}
\eeq

\noi (iii) Let $0\leq k \leq n$ and $\ul{k}\leq l \leq k$. 
Then  $J(n)$ contains 
$(q-2)^l{n\choose l}\left\{ {{n-l}\choose {k-l}}-{{n-l}
\choose{k-l-1}}\right\}$ 
SSJC's starting at rank $k$ and
having offset $l$.
 
\et

Partially order $\Bfq(n)$ by containment. This gives a graded poset 
with rank given by dimension. 
 
\bt \label{mt3}
There exists a SJB $J(n)$ of $V(\Bfq(n))$ such that

\noi (i) The elements of $J(n)$ are orthogonal with respect to
$\langle , \rangle$ (the standard inner product).

\noi (ii) {\em (Singular Values)} Let $0\leq k \leq n/2 $ and let
$(x_k,\ldots ,x_{n-k})$ be any SJC
in $J(n)$ starting at rank $k$ and ending at rank $n-k$. Then
we have, for $k\leq u < n-k$,
\beq \label{mks}
\frac{\pr x_{u+1} \pr}{\pr x_u \pr} & = &
\sqrt{q^k\kq{u+1-k}\kq{n-k-u}}
\eeq
\et

{\textcolor{black} {\bf \section {Complexity of $\Cq(n)$}}}

In this section we prove Theorem \ref{stnb}. The degree of a vertex $(X,f)$
of $\Cq(n)$ is $|X| + (n-|X|)(q-1)$. Define an operator $deg: V(\Bq(n))\rar
V(\Bq(n))$ by
$$ deg((X,f)) = (|X| + (n-|X|)(q-1))(X,f).
$$
We can now write the Laplacian $L: V(\Bq(n))\rar V(\Bq(n))$ of $\Cq(n)$  as
$ L = deg - U - D,
$
where $U,D$ are the up and down operators on $V(\Bq(n))$.

Let $J(n)$ be a rank complete SSJB of $V(\Bq(n))$ satisfying the conditions 
of Theorem \ref{mt2}. Normalize $J(n)$ to get an orthonormal basis
$J'(n)$. Since the vertex degrees are constant on $\Bq(n)_k$ it follows from
the arguments in Section 2 
that the subspace spanned by each SSJC in $J(n)$ is closed under $L$. Using
parts (ii) and (iii) of Theorem \ref{mt2} we can write down the matrix of $L$ in the
basis $J'(n)$.

Let $0\leq k \leq n$ and  $\ul{k} \leq l \leq k$. Let $(x_k,\ldots
,x_{n+l-k})$ be a SSJC in $J(n)$ starting at rank $k$ and having offset $l$.
Set $v_i = \frac{x_{i}}{\pr x_{i} \pr},\;k\leq i \leq n+l-k$. Let $W$ be the
subspace spanned by $\{v_k , \ldots ,v_{n+l-k}\}$. Then $W$ is invariant
under $L$.

Define $M=M(k,n+l-k,n)$ to be the real, symmetric, tridiagonal matrix  of 
size $n+l-2k+1$,
with rows and columns indexed by the set $\{k,k+1,\ldots ,n+l-k\}$, that is
the matrix of $L:W\rar W$ with
respect to the (ordered) basis $\{v_k ,\ldots ,v_{n+l-k}\}$ 
(we take coordinate vectors with respect to a basis as column vectors).
We have from Theorem \ref{mt2} that, for $k\leq i,j \leq n+l-k$, 
the entries of this matrix are given
by: 
\beqn
M(i,j) &=& \left\{ \ba{ll}
      -\sqrt{(q-1)(j-k)(n+l-k-j+1)} & \mbox{if $i=j-1$}\\
                                 & \\ 
      j+ (q-1)(n-j)  & \mbox{if $i=j$}\\
                    & \\
      -\sqrt{(q-1)(j+1-k)(n+l-k-j)} & \mbox{if $i=j+1$}\\
        & \\
      0 & \mbox{if $|i-j|\geq 2$}
                   \ea \right.
\eeqn
It now follows from Theorem \ref{mt2} that the matrix of
$L$ wrt (a suitable ordering of) $J'(n)$ is in block diagonal form,
with
blocks $M(k,n+l-k,n)$, for all $0\leq k \leq n$ and $\ul{k} \leq l \leq k$,
and each such block is repeated 
$(q-2)^l{n\choose l}\left\{ {{n-l}\choose {k-l}}-{{n-l}
\choose{k-l-1}}\right\}$ times. The number of distinct blocks $|\{(k,l) :
0\leq k \leq n,\;\ul{k}\leq l \leq k\}|$ can be easily shown to be 
$(\lfloor \frac{n}{2} \rfloor + 1)( \lceil \frac{n}{2} \rceil + 1)$.  
We now determine the eigenvalues of these blocks. 
In the lemma below the
rows and columns of the matrices on the two sides of the identity are
indexed by different sets (of the same cardinality) but the intended meaning
is clear.

\bl \label{sl} We have
$M(k,n+l-k,n) = (qk - (q-1)l)I + M(0,n+l-2k,n+l-2k).$
\el
\pf From the formula displayed above for the entries of
$M(k,n+l-k,n)$ it follows that the off diagonal entries of the matrices on
both sides of the equation above are the same. The $i$th diagonal entry of
$M(k,n+l-k,n)$ is $k+i-1+(q-1)(n-k-i+1)$ and the $i$th diagonal entry of
$M(0,n+l-2k,n+l-2k)$ is $i-1+(q-1)(n+l-2k-i+1)$ and their difference is $qk -
(q-1)l$, completing the proof. $\Box$

\bt \label{mt4}
The eigenvalues of $M(k,n+l-k,n)$ are $qj - (q-1)l,\;j=k,\ldots
,n+l-k.$
\et
\pf By Lemma \ref{sl} it is enough to show that the eigenvalues of
$M(0,n,n)$ are $qt,\;t=0,\ldots ,n$. 
We shall do this by working with a
suitable linear mapping model for $M(0,n,n)$ and explicitly writing
out the eigenvectors.

Let $J(n)$ be a rank complete SSJB of $V(\Bq(n))$ satisfying the conditions
of Theorem \ref{mt2} and let $J'(n)$ be its normalization. Put
$$
v_k = \sum_{(X,f)\in \Bq(n)_k} (X,f),\;\;v_k' = \frac{v_k}{\pr v_k
\pr},\;\;0\leq k \leq n,
$$
and  set $W=\mbox{ span }\{v_0' ,\ldots ,v_n' \}$. It is easily seen (using
the fact that the bipartite graph between two adjacent ranks of the poset
$\Bq(n)$ is regular on both sides) that the normalization of the unique SSJC
in $J(n)$ starting at rank 0 is $(v_0' ,\ldots ,v_n')$. Thus $W$ is
$L$-invariant and it follows from Theorem \ref{mt2} that the matrix of
$L:W\rar W$ wrt the (ordered) basis $\{v_o' ,\ldots ,v_n'\}$ is $M(0,n,n)$.

Fix $0\leq t \leq n$. Define the vector $w_t \in V(\Bq(n))$ as follows:
\beq
w_t &=& \sum_{(X,f)\in \Bq(n)} \left\{ \sum_{(Y,g)} (-1)^{|X\cap Y|}
\right\} (X,f),
\eeq
where the inner sum is over all $(Y,g)\in \Bq(n)_t$ satisfying: $a\in X\cap
Y$ implies $f(a)=g(a)$.
It is easily seen that $w_t \in W$. We claim that
$L(w_t)=qtw_t$. To prove the claim we introduce a notational device.

The coefficient of $x^k$ in a polynomial $f(x)$ is denoted $[x^k](f(x))$.
The derivative of $f(x)$ is denoted ${\cal D}(f(x))$. We have
$[x^{k-1}]({\cal D}(f(x)))=k([x^k](f(x)))$.

Fix $(X,f)\in \Bq(n)$ with $|X|=k$. Then the coefficient of $(X,f)$ in $w_t$
equals

$$\sum_{j=0}^t (-1)^j {k\choose j}{{n-k}\choose {t-j}} (q-1)^{t-j}
= [x^t]((1-x)^k (1+(q-1)x)^{n-k}).$$

Now the coefficient of $(X,f)$ in $L(w_t)=(deg - U - D)(w_t)$ is equal to
\beqn
&&(k+(n-k)(q-1))\left\{[x^t]((1-x)^k (1+(q-1)x)^{n-k})\right\}\\
&& - (n-k)(q-1)\left\{[x^t]((1-x)^{k+1} (1+(q-1)x)^{n-k-1})\right\}\\
&& - k\left\{[x^t]((1-x)^{k-1} (1+(q-1)x)^{n-k+1})\right\}\\
&&\\
&=&(n-k)(q-1)\left\{[x^t]((1-x)^k (1+(q-1)x)^{n-k-1}
                 (1+(q-1)x -(1-x)))\right\}\\
&& - k \left\{[x^t]((1-x)^{k-1} (1+(q-1)x)^{n-k}
                 (1+(q-1)x -(1-x)))\right\}\\
&&\\
&=&  q(n-k)(q-1)\left\{[x^{t-1}]((1-x)^k (1+(q-1)x)^{n-k-1})\right\}\\
&&  - qk \left\{[x^{t-1}]((1-x)^{k-1} (1+(q-1)x)^{n-k})\right\}\\
&&\\
&=& q \left\{[x^{t-1}]({\cal D}((1-x)^k (1+(q-1)x)^{n-k}))\right\}\\
&&\\
&=& qt \left\{[x^t]((1-x)^k (1+(q-1)x)^{n-k})\right\}
\eeqn
That completes the proof.$\Box$

{\bf Proof} ({\em of Theorem \ref{stnb}}) According to the matrix tree
theorem (see Theorem 5.6.8 in {\bf\cite{s3}}) $c(\Cq(n))$ equals
$\frac{1}{q^n}$ times the product of the nonzero eigenvalues of the
Laplacian of $\Cq(n)$. The graph $\Cq(n)$ being connected the eigenvalue 0
has multiplicity 1 and thus comes from the block $M(0,n,n)$. The product of
the nonzero eigenvalues of $M(0,n,n)$ is $q^n n!$. The result now follows
from  Theorems \ref{mt2}(iii) and \ref{mt4}. $\Box$

{\textcolor{black} {\bf \section {Complexity of $\Cfq(n)$}}}

In this section we prove Theorem \ref{stvs}. The main step of the proof
is the same as that for Theorem \ref{stnb} and thus we will skip some of the
details. For $0\leq k \leq n/2$, define a real, symmetric, tridiagonal matrix
$N=N(k,n-k,n)$ of size $n-2k+1$,
with rows and columns indexed by the set $\{k,k+1,\ldots ,n-k\}$, and with
entries given as follows. 

For $k\leq i,j \leq n-k$ define
\beqn
N(i,j) &=& \left\{ \ba{ll}
      -\sqrt{q^k\kq{j-k}\kq{n-k-j+1}} & \mbox{if $i=j-1$}\\
                                 & \\ 
      \kq{j} + \kq{n-j}  & \mbox{if $i=j$}\\
                    & \\
      -\sqrt{q^k\kq{j+1-k}\kq{n-k-j}} & \mbox{if $i=j+1$}\\
        & \\
      0 & \mbox{if $|i-j|\geq 2$}
                   \ea \right.
\eeqn
For $0\leq k \leq n/2 $ and $k\leq j \leq n-k+1$ define
$N_j=N_j(k,n-k,n)$
to be the principal submatrix of $N=N(k,n-k,n)$ indexed by the rows and
columns in the set $\{j,j+1,\ldots ,n-k\}$. Thus, $N_k = N$ and $N_{n-k+1}$
is the empty matrix, which by convention has determinant 1.

\bl \label{det}
 
For $0\leq k \leq n/2 $ and $k\leq j \leq n-k+1$ we have 
 
(i) $F_q(n,k,j) = \mbox{\em det}(N_j(k,n-k,n)).$

(ii) $F_q(n,0,j) = \kq{n}\kq{n-1}\cdots \kq{j}.$

(iii) The degree of $F_q(n,0,j)$ is $\sum_{t=j}^n (t-1)$
and, for $k\geq 1$, the degree of $F_q(n,k,j) = \sum_{t=j}^{n-k}\mbox{\em
max}\{t-1,n-t-1\}$.
\el
\pf (i) By (reverse) induction on $j$. The base cases $j=n-k+1,n-k$ are 
clear and the general case follows by expanding the determinant of $N_j$
along its first column.

(ii) By (reverse) induction on $j$. The base cases $j=n+1,n$ are clear. By
induction and the defining recurrence for $F_q(n,k,j)$ we have 
\beqn
F_q(n,0,j) & = & (\kq{j} + \kq{n-j})F_q(n,0,j+1) -
(\kq{j+1}\kq{n-j})F_q(n,0,j+2)\\
& = & (\kq{j} + \kq{n-j})\kq{n}\cdots \kq{j+1} -
(\kq{j+1}\kq{n-j})\kq{n}\cdots \kq{j+2}\\
&=& \kq{n}\cdots \kq{j}. 
\eeqn

(iii) The degree of $F_q(n,0,j)$ follows from part (ii) above. Now assume
that $k\geq 1$. We prove the stated formula by (reverse) induction on $j$.
The formula clearly holds for $j=n-k+1,n-k$. By the inductive hypothesis the
degree of the first term on the rhs of the defining recurrence (\ref{dr}) is
$\sum_{t=j}^{n-k} \mbox{ max}\{t-1,n-t-1\}$ and the 
degree of the second term on the rhs of (\ref{dr}) is
$n-k-1 + \sum_{t=j+2}^{n-k} \mbox{ max}\{t-1,n-t-1\}$. The result will be
proven if we show that
$r(n,j)= \mbox{ max}\{j-1,n-j-1\} + \mbox{ max}\{j,n-j-2\} > n-k-1.$
But this is clear, since $\mbox{max}\{t-1,n-t-1\} \geq n/2 -1$ and thus
$r(n,j)\geq n/2 - 1 + n/2 = n-1 > n-k-1$, since $k>1$.$\Box$

{\bf Proof} {\em (of Theorem \ref{stvs})} Let $J(n)$ be a SJB of
$V(\Bfq(n))$ satisfying the conditions of Theorem \ref{mt3}. Normalize
$J(n)$ to get an orthonormal basis $J'(n)$. Let $L$ denote the Laplacian of
$\Cfq(n)$. Just as in the case of $\Cq(n)$ in Section 3, it follows from
Theorem \ref{mt3} that the matrix of
$L$ wrt (a suitable ordering of) $J'(n)$ is in block diagonal form,
with
blocks $N(k,n-k,n)$, for all $0\leq k \leq n/2$ 
and each such block is repeated 
$\qb{n}{k} - \qb{n}{k-1}$ times. The number of distinct blocks is $1 +
\lfloor n/2 \rfloor$.

The unique element in $J'(n)$ of rank 0 is the vector ${\bf 0}$ (here ${\bf
0}$ is the zero subspace).

Let ${\cal M}$ denote the matrix of the Laplacian of $\Cfq(n)$ in the standard
basis $\Bfq(n)$ and let ${\cal M}'$ be obtained from ${\cal M}$ by removing the row and column
corresponding to vertex ${\bf 0}$.
According to the matrix tree
theorem (see Theorem 5.6.8 in {\bf\cite{s3}}) 
$c(\Cfq(n))=\mbox{ det}({\cal M}')$.
A little reflection shows that, by changing bases from 
$\Bfq(n) - \{{\bf 0}\}$ to $J'(n) - \{{\bf 0}\}$, 
${\cal M}'$ block diagonalizes
with a block $N_1(0,n,n)$ of multiplicity 1 and blocks $N(k,n-k,n)$, for all
$1\leq k \leq n/2$, of multiplicity $\qb{n}{k} - \qb{n}{k-1}$. The result
now follows from Lemma \ref{det}.$\Box$ 

{\bf Remark} A natural question at this point is whether there is a vector
space analog of Theorem \ref{mt4}. This would involve guessing the
(eigenvalue, eigenvector) pairs of $N(k,n-k,n)$ and then verifying this guess
using an analog of the proof of Theorem \ref{mt4}.

{\textcolor{black} {\bf \section {Orthogonal SSJB  
of $\Bq (n)$ and SJB of $\Bfq(n)$}}}

In this section we prove the results stated in Section 2. We begin with the
proof of Theorem \ref{mt3}.

{\bf Proof} ({\em of Theorem \ref{mt3}}) We shall put together several 
standard results.
 
(i) The map $U^{n-2k}:V(\Bfq(n)_k)\rar
V(\Bfq(n)_{n-k}),\;0\leq k \leq n/2$ is well known to be bijective. 
It follows, using a standard Jordan canonical
form argument, that an SJB of $V(\Bfq(n))$ exists.

(ii) Now we show existence of an orthogonal SJB. We use the action 
of the group $GL(n,\Fq)$ on $\Bfq(n)$.
As is easily seen the
existence of an orthogonal SJB of $V(\Bfq(n))$ (under the standard inner
product)
follows from facts (a)-(d) below  by an application of Schur's lemma:

\noi (a) Existence of some SJB of $V(\Bfq(n))$.

\noi (b) $U$ is $GL(n,\Fq)$-linear.

\noi (c) For $0\leq k \leq n$, $V(\Bfq(n)_k)$ is the sum of
$\mbox{min}\{k,n-k\}+1$ distinct
irreducible $GL(n,\Fq)$-modules (this result is well known. The
corresponding result for the $S_n$ action on $V(B(n)_k)$ is proved in
Chapter 29 of {\bf\cite{jl}}. An identical proof works in the present case).

\noi (d) For a finite group $G$, a $G$-invariant inner product on an
irreducible $G$-module is unique upto scalars.

(iii) Now we prove part (ii) of Theorem \ref{mt3}. 
Define an operator $H: V(\Bfq(n))\rar V(\Bfq(n))$ by
$$H(X)= (\kq{k} - \kq{n-k})X,\;\; X\in \Bfq(n)_k,\;0\leq k \leq n.$$ 
It is easy to check that $[U,D] = UD - DU = H$. To see this, fix $X\in
\Bfq(n)_k$, and note that 
$UD(X) = \kq{k} X + \sum_Y Y$, where the sum is
over all $Y\in \Bfq(n)_k$ with $\mbox{dim}(X\cap Y) = k-1$. 
Similarly, $DU(X) = \kq{n-k} X + \sum_Y Y$, where the sum is
over all $Y\in \Bfq(n)_k$ with $\mbox{dim}(X\cap Y) = k-1$. 
Subtracting we get  $[U,D]=H$.

Let $J(n)$ be an orthogonal SJB of $V(\Bfq(n))$ and let $(x_k,\ldots
,x_{n-k})$ be a SJC in $J(n)$ starting at rank $k$ and ending at rank $n-k$.
Put
$x_j' = \frac{x_j}{\pr x_j \pr}$ and $\alpha_j = \frac{\pr x_{j+1} \pr}{\pr
x_{j} \pr},\;k\leq j \leq n-k$ (we take $x_{k-1}'=x_{n-k+1}'=0$). We have, 
from (\ref{trick}) and (\ref{trick1}),
$$U(x_j')=\alpha_j x_{j+1}',\;\;D(x_{j+1}')=\alpha_j x_j',\;\;k\leq j <
n-k.$$

We need to show that 
\beq \label{ind}
\alpha_j^2 &=& q^k \kq{j+1-k}\kq{n-k-j},\;\;k\leq j < n-k .
\eeq 

We show this by
induction on $j$. We have $DU = UD - H$. Now $DU(x_k')=\alpha_k D(x_{k+1}')
= \alpha_k^2 x_k'$ and $(UD-H)(x_k')= (\kq{n-k} - \kq{k}) x_k'$ (since
$D(x_k')=0$). Hence
$\alpha_k^2 = \kq{n-k} - \kq{k} = q^k \kq{n-2k}$. Thus (\ref{ind}) holds for
$j=k$.

As in the previous paragraph $DU(x_j')=\alpha_j^2 x_j'$ and
$(UD-H)(x_j')=(\alpha_{j-1}^2 + \kq{n-j} - \kq{j}) x_j'$.
By induction, we may assume $\alpha_{j-1}^2 = q^k \kq{j-k}\kq{n-k-j+1}$.
Thus we see that $\alpha_j^2$ is   
\beqn
&=& q^k \kq{j-k}\kq{n-k-j+1} + \kq{n-j} - \kq{j}\\
&=& q^k \left\{ (\kq{j+1-k} - q^{j-k})(\kq{n-k-j} + q^{n-k-j})\right\} 
+ \kq{n-j} - \kq{j}\\
&=& q^k \left\{ \kq{j+1-k}\kq{n-k-j} + q^{n-k-j}\kq{j+1-k} 
 -q^{j-k}\kq{n-k-j} - q^{n-2k}\right\} \\
&& + \kq{n-j} - \kq{j}\\
&=& q^k \kq{j+1-k}\kq{n-k-j} + q^{n-j}\kq{j+1-k} - q^j\kq{n-k-j} - q^{n-k}\\
&&     + \kq{n-j} - \kq{j}\\
&=& q^k \kq{j+1-k}\kq{n-k-j} + \kq{n+1-k} - \kq{n-j} - \kq{n+1-k} + \kq{j}\\
  &&   + \kq{n-j} - \kq{j}\\
&=& q^k \kq{j+1-k}\kq{n-k-j},
\eeqn
completing the proof.$\Box$

\noi {\bf Remark} The $q=1$ case of Theorem \ref{mt3} yields Theorem
\ref{mt1}. In {\bf\cite{s1}} a constructive proof of Theorem \ref{mt1} was
given by producing an explicit orthogonal SJB of $V(B(n))$, together with a
representation theoretic interpretation of this basis. It would be 
interesting to construct an explicit orthogonal SJB of $V(\Bq(n))$.

Now we prove Theorem \ref{mt2}. Consider the following identity
\beq \label{ubd}
q^n &=& (q-2+2)^n \,\;=\,\; \sum_{l=0}^n {n\choose l} (q-2)^l 2^{n-l}.
\eeq
We shall give a linear algebraic interpretation to the identity above, 
which reduces Theorem \ref{mt2} to Theorem \ref{mt1}. We
begin with a combinatorial interpretation of (\ref{ubd}) which suggests the
algebraic interpretation.

A subset $S\subseteq \Bq(n)$ is said to be {\em upper Boolean of rank $t$}
if $\mbox{rankset}(S)=\{n-t,n-t+1,\ldots ,n\}$ and $S$, 
with the induced order, is order isomorphic to a Boolean
algebra $B(t)$. 

\bt There is a partition of $\Bq(n)$ into pairwise disjoint upper
Boolean subsets, with $(q-2)^l {n\choose l}$ of them  having rank $n-l$, for
each $l=0,1,\ldots ,n$.
\et
\pf Let $0\leq l \leq n$ and let $X\subseteq [n]$ with $|X|=l$. Fix $f:X
\rar \{1,\ldots ,q-1\}$ with $f(a) \not= 1$ for all $a\in X$. Let
$\Bq(n,l,X,f)$ denote the set of all $(Y,g)\in \Bq(n)$ with $X\subseteq Y$,
$f(a) = g(a),\;a\in X$, and $g(a)=1,\;a\in Y-X$. Clearly $\Bq(n,l,X,f)$ is
an upper Boolean subset of rank $n-l$. Once $l$ is fixed, $X$ can be chosen
in ${n\choose l}$ ways and then $f$ can be chosen in $(q-2)^l$ ways. Going
over all choices of $l,X,f$ we get the required decomposition.$\Box$ 

Let $(V,f)$ be a pair consisting of a finite dimensional inner product space
$V$ (over $\C$) and a linear operator $f$ on $V$. Let $(W,g)$ be
another such pair. By an isomorphism of pairs $(V,f)$ and $(W,g)$ we mean a
linear isometry (i.e, an inner product preserving isomorphism) $\theta : V
\rar W$ such that $\theta(f(v)) = g(\theta(v)),\;v\in V$.

Consider the inner product space $V(\Bq(n))$, with the standard inner
product. An {\em upper Boolean subspace} of rank $t$ is a homogeneous
subspace $W\subseteq V(\Bq(n))$ such that
$\mbox{rankset}(W)=\{n-t,n-t+1,\ldots ,n\}$, $W$ is closed under the up
operator $U$, and there is an isomorphism of pairs $(V(B(t)),\sqrt{q-1}\,U)
\cong (W,U)$ that sends homogeneous elements to homogeneous elements and
increases rank by $n-t$ (here, and in the rest of this section,  we use $U$ 
to denote the up operator on both $V(B(n))$ and $V(\Bq(n))$. The context
always makes clear which poset is intended).

\bt \label{oubd} There exists an orthogonal decomposition of $V(\Bq(n))$ 
into upper Boolean subspaces, with $(q-2)^l {n\choose l}$ of 
them  having rank $n-l$, for each $l=0,1,\ldots ,n$.
\et

Before proving Theorem \ref{oubd} let us see how it implies Theorem
\ref{mt2}.

{\bf Proof} ({\em of Theorem \ref{mt2}}) Take an orthogonal
decomposition of $V(\Bq(n))$ into upper Boolean subspaces given by Theorem
\ref{oubd} and let $W$ be an
upper Boolean subspace in this decomposition of rank $n-l$.

Use Theorem \ref{mt1} to get an orthogonal SJB $J(n-l)$ of $V(B(n-l))$ 
wrt $\sqrt{q-1}\;U$ 
(rather than just $U$) and transfer it to $W$. 
Each SJC in $J(n-l)$ will get transfered
to a SSJC in $V(\Bq(n))$ of offset $l$ and, using (\ref{mks1}), we
see that this SSJC will satisfy (\ref{mks2}). The number of these SSJC's (in
$W$)
starting at rank $k$ is ${{n-l}\choose {k-l}} - {{n-l}\choose {k-l-1}}$ and
since the number of rank $n-l$ upper Boolean subspaces in the decomposition
is $(q-2)^l{n \choose l}$, Theorem \ref{mt2}  now follows.$\Box$ 

Fix a $(q-1) \times (q-1)$ unitary matrix $A=(a_{ij})$, with rows and columns indexed by
$\{1,2,\ldots ,q-1\}$, and with first row $\frac{1}{\sqrt{q-1}}(1,1,\ldots
,1)$. 

We now prove Theorem \ref{oubd}.

{\bf Proof} ({\em of Theorem \ref{oubd}}) We give an inductive procedure to
explicitly construct an orthogonal decomposition of $V(\Bq(n))$ into upper
Boolean subspaces. The case $n=0$ is clear. In this proof it is convenient
to think of the elements of $\Bq(n)$ as $n$-tuples $(x_1,\ldots ,x_n)$ with
$x_i \in \{0,1,\ldots ,q-1\}$ for all $i$. 

Consider $V=V(\Bq(n+1))$ with the standard inner product. 
Define $W(0)$ to be the subspace of $V$ spanned by
all elements $(x_1,\ldots,x_n,0)\in \Bq(n+1)$ with last coordinate $0$ and
define $V'$ to be the subspace of $V$ spanned by
all elements $(x_1,\ldots,x_n,i)\in \Bq(n+1)$ with 
$i\in \{1,2,\ldots ,q-1\}$. We have an orthogonal direct sum decomposition
$V=W(0) \oplus V'$.

For $i=1,2,\ldots ,q-1$ define linear maps
$$ L_i : W(0) \rar V' $$
by $L_i((x_1,\ldots ,x_n,0)) = \sum_{j=1}^{q-1}a_{ij}(x_1,\ldots ,x_n,j).$
Set $W(i) = \mbox{ Image}(L_i),\;1\leq i \leq q-1$. It is easy to see from
the definition of $A$ that

(i) $L_i: W(0) \rar W(i)$ is an isometry, for all $i$.

(ii) $V=W(0)\oplus V' = W(0)\oplus W(1) \oplus W(2) \oplus \cdots \oplus
W(q-1)$ is an orthogonal direct sum decomposition.

There is an isometry $V(\Bq(n))\cong W(0)$ given by $(x_1,\ldots
,x_n)\mapsto (x_1,\ldots ,x_n,0)$. We denote the up operator on $V(\Bq(n))$
by $U_n$ and the corresponding operator on $W(0)$ under the isometry above
by the same symbol $U_n$. We denote the up operator on $V(\Bq(n+1))$ by
$U_{n+1}$. A little reflection shows that
\beq \label{cl}
U_{n+1}(L_i(x_1,\ldots ,x_n,0))&=&L_i(U_n(x_1,\ldots ,x_n,0)),\;1\leq i \leq
q-1, \\ \label{upi}
U_{n+1}((x_1,\ldots ,x_n,0))&=& U_n((x_1,\ldots ,x_n,0)) +
\sqrt{q-1}\;L_1((x_1,\ldots ,x_n,0)).
\eeq
It follows from (\ref{cl}) above that $W(i)$ is closed under $U_{n+1}$,
for $1\leq i \leq q-1$.

As in the paragraph above we identify the pair $(W(0),U_n)$ with
$(V(\Bq(n)),U_n)$.
Let $X\subseteq W(0)$ be an upper Boolean subspace of rank
$t$. 
By (\ref{cl}) above, there is an
isomorphism of pairs $(X,U_n)\cong (L_i(X),U_{n+1}),\;1\leq i \leq q-1$.
Since $L_i$ increases rank by 1, it follows that each $L_i(X),\;1\leq i \leq
q-1$, is an upper Boolean subspace of $V(\Bq(n+1))$ of rank $t$.

Now we claim that $X\oplus L_1(X)$ is an upper Boolean subspace of
$V(\Bq(n+1))$ of rank $t+1$. We show that the pair $(X\oplus
L_1(X),U_{n+1})$ has the same recursive structure as the pair
$(V(B(t+1)),\sqrt{q-1}\;U_{t+1})$. 

Define $R_1: V(B(t)) \rar V(B(t+1))$ by $R_1(s)=s\cup \{t+1\},\;s\in B(t)$.
We have 

(a) $V(B(t+1)) = V(B(t)) \oplus R_1(V(B(t)))$ is an orthogonal
decomposition.

(b) $R_1$ is an isometry.

(c) $R_1(V(B(t)))$ is closed under $U_{t+1}$ and $R_1:V(B(t)) \rar
R_1(V(B(t)))$ is an isomorphism of pairs $(V(B(t)),\sqrt{q-1}\;U_t) \cong
(R_1(V(B(t))),\sqrt{q-1}\;U_{t+1})$.

(d) $\sqrt{q-1}\;U_{t+1}(s) = \sqrt{q-1}\;U_t(s) + \sqrt{q-1}\;R_1(s),\;s\in
B(t)$.

The corresponding statements about $L_1$ are

(a') $X \oplus L_1(X)$ is an orthogonal
decomposition.

(b') $L_1$ is an isometry.

(c') $L_1(X)$ is closed under $U_{n+1}$ and $L_1:X \rar
L_1(X)$ is an isomorphism of pairs $(X,U_n) \cong
(L_1(X),U_{n+1})$.

(d') $U_{n+1}(v) = U_n(v) + \sqrt{q-1}\;L_1(v),\;v\in X$.

The claim easily follows from statements (a)-(d) and (a')-(d') above.

So, from an upper Boolean subspace $X$ of rank $t$ in $W(0) (\cong
V(\Bq(n)))$ we get one
upper Boolean subspace $X\oplus L_1(X)$ of rank $t+1$ and $q-2$ upper
Boolean subspaces $L_i(X),\;i=2,\ldots ,q-1$ of rank $t$ in $V(\Bq(n+1))$.

Now, using the inductive hypothesis take an 
orthogonal decomposition of $V(\Bq(n))$ 
into upper Boolean subspaces, with $(q-2)^l {n\choose l}$ of 
them  having rank $n-l$, for each $l=0,1,\ldots ,n$. For each upper Boolean
subspace in this decomposition produce upper Boolean subspaces in
$V(\Bq(n+1))$ as in the paragraph above. Clearly, this will give a
orthogonal decomposition of $V(\Bq(n+1))$. The number of upper Boolean
subspaces of rank $n+1-l$ in this decomposition is
$$(q-2)(q-2)^{l-1}{n\choose {l-1}} + (q-2)^{l}{n\choose {l}} = (q-2)^l
{{n+1}\choose l},$$
completing the proof. $\Box$

\begin{center} {\bf \Large{Acknowledgements}}
\end{center}
It is a pleasure to thank Navin Singhi for several stimulating
discussions concerning the topic of this paper and the papers
{\bf\cite{s1,s2}}. I thank Thomas Zaslavsky for several useful suggestions
on the write up. I am grateful to Professor Alexander Schrijver for his
encouragement.

\end{document}